\pdfoutput=1
\documentclass{article}

\usepackage[T1]{fontenc}
\usepackage[UKenglish]{babel}

\usepackage{nicefrac, amsmath, amsthm, amssymb, mathabx}
\usepackage{Factory, Classic, Math, Styl, Theorema}

\usepackage[all, 2cell]{xy} \UseAllTwocells 
\newcommand{\upar}{ \ar@<+.4ex> }			\newcommand{\downar}{ \ar@<-.4ex> }
\newcommand{\curvear}{\ar@/^1.5pc/} 		\newcommand{\cruvear}{\ar@/^-1.5pc/}	
\xymatrixcolsep{4pc} 

\newcommand{\K}{\mathbf{K}}

\newcommand{\tbeta}{\tilde{\beta}}
\newcommand{\AMod}{{}_A\mathfrak{Mod}}

\usepackage{fancyhdr}
\begin{document}

\lefthyphenmin=2 \righthyphenmin=2

\titul{GABRIEL 2-QUIVERS FOR \\ FINITARY 2-CATEGORIES}
\auctor{Qimh Richey Xantcha%
\footnote{\textsc{Qimh Richey Xantcha} (\texttt{qimh@math.su.se}), Stockholm University.}}
\datum{\today}
\maketitle

\bigskip 

\epigraph{
Fru Lenngren har hedrat ålderdomen såsom sådan, hon har skrifvit tvenne vackra kväden [\dots] 
till den hedersamma ålderdomens pris, men den ålderdom, som skämdes för sin ålder, mot den 
har hon ett helt koger af pilar [\dots]. 
\attr{Från Warburgs Företal til Lenngrens \emph{Skaldeförsök}}
}

\bigskip 

\begin{argument}
\noindent We develop the theory of 2-quivers and quiver 2-categories to run in parallel 
with the classical theory of quiver algebras. A quiver 2-category is always finitary, 
and, conversely, 
every finitary 2-category will be bi-equivalent with a quiver 2-category for a unique underlying reduced 2-quiver. 

As an application, we produce an example of a fiat 2-category, one of whose Duflo involutions is not 
self-adjoint.   

\end{argument}

\bigskip

\noindent
Gabriel's combinatorial approach to the study of 
finite\hyp dimensional algebras via quivers
(see, for example, \cite{Gabriel}), plays a central rôle in modern representation theory. Its success is due to 
the following facts about finite\hyp dimensional, associative algebras 
over an algebraically closed field $\K$: 
 
\bnum
\item Every finite-dimensional algebra is Morita equivalent to a basic algebra, unique up to isomorphism 
(where \emph{basic} means that all simple modules are 1\hyp dimensional).
\item Every \emph{quiver algebra} --- a quotient of a path algebra by an admissible ideal --- 
is finite-dimensional and basic.  
\item Conversely, every finite-dimensional, basic algebra is isomorphic with a quiver algebra, 
for a unique underlying quiver, known as the \emph{Gabriel quiver} of the algebra. 
\enum

Thus, from the point of view of representation theory, there is no loss of generality in considering 
quiver algebras only. 
One major advantage is that quiver algebras are well suited 
to various computational 
purposes, in particular, for constructing examples.

Recent years have seen a progressive trend towards $\K$-linear 2-categories and their representations. 
One of the main reasons is that the procedure of \emph{categorification} 
is now usually clad in the formalism of 2-categories. The categorification approach was recently 
used  to prove various spectacular results in different areas of mathematics. 
Here one could mention 
Khovanov's categorification of the Jones polynomial, see \cite{Khovanov}; Chuang \& Rouquier's 
proof of Broué's conjecture for blocks of the symmetric group,  see \cite{Chuang}; 
the categorification of the Temperley--Lieb algebra in \cite{Bernstein}, \cite{Stroppel} and many others. 
For a systematic account of algebraic categorification and its applications we refer to Mazorchuk's treatise \cite{Categorification} on the subject.
The $2$-representation theory of finitary $2$-categories is studied in a series \cite{Cell}, \cite{Abelian}, \cite{Endo}, \cite{Morita} of papers by Mazorchuk \& Miemietz.

Since a finite-dimensional algebra is essentially a \emph{finitary} 1-category 
(see Section \ref{A: Prel}), 
a finitary 2-category would be a natural analogue at the next level. 
We refer to Section~\ref{A: Prel} for a brief 
review of 1-representation theory from a categorical viewpoint.

It would then seem natural to request a corresponding theory of \emph{2-quivers}, capable of encoding 
finitary 2-categories. This is the purpose of the present paper.

We show that, under the apposite modifications, 
Gabriel's classical results indeed generalise. To this end, we define 
(finite) 2-quivers below (Definition \ref{D: 2-quiver}) in the ``obvious'' way. 
The construction of quiver 2-categories then proceeds in four steps:
\bnum
\item Factor out relations on the 1-level to form a \emph{bound 2-quiver}. (Definition~\ref{D: Bound}) 
\item Create a \emph{path 1\nicefrac12-category} based upon this bound 2-quiver; 
it is like a 2-category, except there is no horizontal composition. (Definition~\ref{D: Path Cat})
\item Factor out an \emph{admissible} ideal to form a quotient 1\nicefrac12-category.  
(Definition \ref{D: Admissible})
\item Define horizontal composition to create a 2-category. 
This is accomplished through \emph{vertical drop relations}. (Definition \ref{D: Drop})
\enum

We may now state our principal results. 

\begin{inttheorem}[\ref{S: Finitary}]
Quiver 2-categories are finitary, the indecomposable 1-cells being precisely the 1-edges of the underlying 2-quiver.
\end{inttheorem}

\begin{inttheorem}[\ref{S: Every}]
Every finitary 2-category is bi-equivalent with a quiver 2-category for a unique underlying reduced 2-quiver. 
\end{inttheorem}

The last section presents some examples and applications. Worthy of note 
is Example \ref{X: Duflo} 
(p.~\pageref{X: Duflo}), 
which exhibits the first known example of a fiat category with a Duflo involution that is not self-adjoint, 
answering a question posed by Mazorchuk at the 2-seminar at Uppsala.

We are grateful to Volodymyr Mazorchuk for fruitful scientific discussions, as well as his most valuable comments on the draft of this manuscript.

\setcounter{section}{-1} 
\section{1-Representation Theory}			\label{A: Prel}

Throughout, $\K$ denotes a fixed field, always assumed algebraically closed. 
All tensor products shall be taken over this $\K$, and all algebras shall be assumed finite-dimensional, associative, 
and unital $\K$-algebras. All 1- and 2-categories will be $\K$-linear, i.e.~enriched over the category of 
$\K$-vector spaces. We will say that a 
$\K$-linear category is \textbf{finitary}
if it has finitely many objects and all Hom-spaces are finite-dimensional over $\K$.

Studying representations of algebras (associative, unital and finite\hyp dimensional) is equivalent to studying representations 
of finitary categories.
We devote this preliminary section 
to a brief recapitulation of 1-representation theory from the categorical point of view. 
A more comprehensive treatment can be found in \cite{Roiter}. See also Section 2 of \cite{Morita}.

Given a finitary category $C$ with objects $X_1,\dots,X_n$, the space 
$$
\bigoplus_{i,j=1}^n C(X_i,X_j)
$$ 
has the canonical structure
of a $\K$-algebra induced by the composition in $C$ (with any
undefined compositions postulated to be $0$). Conversely, any 
finite-dimensional algebra $A$ defines a 
single-object finitary category in the obvious way, with elements of the algebra being the arrows of the category. 
However, this naïve approach can be slightly refined.

Associate to $A$ a category $C_A$ by choosing as its objects a complete system $e_1,\dots,e_k$
of primitive, orthogonal idempotents  in $A$, so that $1=e_1+\dots+e_k$. Some of these idempotents may
be equivalent, but note that we opt not to discard multiple copies. Set $C_A(e_i,e_j)=e_iAe_j$ and let composition in $C_A$ be induced from the 
multiplication in $A$. This establishes 
a one-to-one correspondence between isomorphism classes of finite-dimensional algebras and finitary categories $C$, 
each of whose endomorphism rings $C(X,X)$ is \emph{local}.
 A dictionary for translating between the two settings is given in Table~\ref{T: Dictionary}.
\btab{cc}
\hline
\textbf{Categories $C_A$} & \textbf{Algebras $A$} \\
\hline
finitary & finite-dimensional \\
connected & indecomposable \\
disjoint union & direct sum \\
skeletal & basic \\
skeletally discrete
& semi-simple \\
skeletally discrete on a single object & simple \\
isomorphic & isomorphic \\
equivalent & Morita equivalent  \\
\hline
\etabc{Dictionary for translating between algebras and their associated categories over an algebraically closed field.
\label{T: Dictionary}} 

Gabriel's fundamental insight was that a finite-dimensional algebra, or, equivalently, a finitary category, 
may be compressed into a finite collection of data, viz.~a (finite) quiver with (admissible) relations. 
This is the \emph{Gabriel quiver} of the algebra. Details on its construction, and the recipe for 
re-assembling the original category from the combinatorial data, are described, for instance, in 
\cite{ARS}. Stated in categorical language, Gabriel's results assert that 
quiver categories are finitary and, conversely, every finitary category over an algebraically closed field $\K$ 
is isomorphic with a quiver category for a unique underlying quiver. 

Let us recall the well-known properties of the Jacobson radical of a finitary category $C$, which will be frequently used in the sequel: 
\bnum
\item If $R$ is a nilpotent ideal of $C$ such that the quotient category $C/R$ is discrete, then $R=\Rad C$.
(Conversely, the radical is always nilpotent, but, unless $\K$ be algebraically closed, 
$C/\Rad C$ may not be discrete, a counter-example being given by the 
$\R$-algebra $\C$.)

\item
Any radical element can be expressed as a polynomial 
in the vectors of an arbitrary basis of $\Rad C/\Rad^2 C$. 
\enum

\section{2-Categories}

All our 2-categories shall be strict 
and $\K$-linear. The informal definition we give below is intended to serve as a reminder, but 
also to introduce the typographical 
conventions which will be strictly observed throughout the paper. 
For more details on 2-categories, the reader is referred to Leinster's note \cite{Leinster}.

\bdf
A \textbf{linear 2-category} is a category enriched over the category of linear categories. 
It consists of: 
\begin{description}
\item[0-cells.] Denoted by the letters $X$, $Y$, $Z$, \dots 
\item[1-cells.] Denoted by $F$, $G$, $H$, \dots
\item[2-cells.] Denoted by $\alpha$, $\beta$, $\gamma$, \dots~(forming vector spaces). 
\item[1-cell composition.] Denoted by $FG=F\circ G$.
\item[Horizontal 2-cell composition.] Denoted by $\alpha\beta=\alpha\circ_0\beta$.
\item[Vertical 2-cell composition.] Denoted by $\alpha\circ\beta=\alpha\circ_1\beta$.  
\item[Identity 1-cells.] Denoted by $I_X$, $I_Y$, $I_Z$, \dots 
\item[Identity 2-cells.] Denoted by 
$\iota_{I_X}$, $\iota_{I_Y}$, $\iota_{I_Z}$, $\iota_F$, $\iota_G$, $\iota_H$, \dots
\end{description}
All three compositions are to be \emph{(strictly) associative}, that on the 
1-level \emph{biadditive} and those on the 2-level \emph{bilinear}. Moreover, 
horizontal composition is required to be bifunctorial, which, in concrete terms, means 
it should satisfy the equality $\iota_F\iota_G=\iota_{FG}$ (assuming compatibility)
as well as the \emph{Interchange Law}
$$
(\alpha \circ \beta) (\alpha'\circ \beta') = (\alpha\alpha')\circ (\beta\beta')
$$ 
for merging 2-cells: 
$$ \xymatrixcolsep{4pc} 
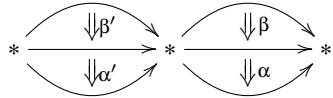
\xymatrix{
\ast \ruppertwocell{\beta'} \ar[r] \rlowertwocell{\alpha'}
& \ast \ruppertwocell{\beta} \rlowertwocell{\alpha} 
\ar[r] & \ast
}
$$
\edf

\bex
The primordial example of a 2-category is the 2-category of small linear categories. The 1-cells are functors, 
and the 2-cells are natural transformations under their vertical and horizontal composition. 
\eex

Let $C$ be a 2-category. We agree to write $C[n]$ ($n=0,1,2$) for the set of $n$-cells, 
and furthermore $C[0,1]$ for the 1-category obtained by excising the 2-cells of $C$. 
The symbol $C(X,Y)$ denotes the set of 1-cells 
from $X$ to $Y$. Likewise, $C(F,G)$ stands for the linear space of 2-cells from
$F$ to $G$.  

The notion of a finitary 2-category originated with Mazorchuk \& Miemietz in \cite{Cell} (paragraph 2.2). 
It may justly be thought of 
as the 2-level analogue of a finite-dimensional algebra. 

\begin{definition}		\label{D: Finitary}
The 2-category $C$ is \textbf{finitary} if it has: 
\bnum
\item finitely many 0-cells; 
\item fully additive 1-categories $C(X,Y)$
(\emph{id est}: finite direct sums exist and all vertical idempotents split); 
\item finitely many indecomposable 1-cells up to isomorphism; 
\item indecomposable identity 1-cells; 
\item finite-dimensional spaces $C(F,G)$.
\enum 
\end{definition}

Conditions 3 and 5 may be summarised by the requirement that each $C(X,Y)$ be a finitary 1-category. 

Full additivity implies that finitary 2-categories 
possess the Krull--Schmidt property: \emph{Every 1-cell has a  unique decomposition into indecomposables.}
The reader is referred to Appendix A of \cite{Burban2} for details on this. 

Let $X$ and $Y$ be 0-cells of $C$. Then $C(X,Y)$ is a finitary 1-category, and so will possess 
a radical $\Rad C(X,Y)$ with the usual properties. 

\bdf 
The \textbf{(vertical) Jacobson radical} $\Rad C$ is the substructure of $C$ 
having the same 0- and 1-cells as $C$, and 2-cells 
$$
(\Rad C)(F,G) = \Rad C(X,Y) \cap C(F,G) \qquad\text{for $F,G\colon X\to Y$.} 
$$
\edf

$\Rad C$ is a \textbf{vertical ideal} in the sense that 
$\alpha\circ\rho, \rho\circ\beta\in \Rad C$ for any $\alpha,\beta\in C$ 
and $\rho\in\Rad C$ (assuming both compositions defined). It is \emph{not} an horizontal ideal of $C$.

Two 0-cells $X$ and $Y$ of a 2-category $C$ are called \emph{equivalent} if there exist
1-cells $F:X\to Y$ and $G:Y\to X$ such that 
$GF$ is isomorphic to $I_X$ in $C(X,X)$ and 
$FG$ is isomorphic to $I_Y$ in $C(Y,Y)$. 
A functor $\Theta:C\to D$
between two 2-categories $C$ and $D$ is called a \textbf{bi-equivalence} if it is  
surjective up to equivalence on objects and is a \emph{local equivalence}; that is, it induces an
equivalence $C(X,Y)\cong D(\Theta(X),\Theta(Y))$ for all $X$ and $Y$.
Bi-equivalence is the natural notion of ``equivalence'' for 2-categories, see e.g.~\cite{Leinster} or \cite{Morita}.

\section{Licit Ideals}	\label{A: Licit}

A \emph{semi-ring} (also known as a \emph{rig}) is a binary structure satisfying all ring axioms except 
possibly the existence of additive inverses. On the other hand, it is necessary to append an axiom positing that 
$0\cdot a=a\cdot 0=0$ for all $a$. (This law is a theorem for rings, but required as an axiom for semi-rings.)

A major difference between rings and semi-rings becomes apparent in the discussion of factor semi-rings. 
In the general case, these are constructed by dividing out by a \emph{congruence relation}, 
whose equivalence classes are not necessarily cosets with respect to some ideal. 

The situation simplifies for semi-rings in which additive cancellation is allowed.  
A semi-ring $R$ of this type may be embedded in a ring $R^\Delta$, known as the \emph{ring of differences}, 
 via the Grothendieck construction (cf.~Chapter 8 of \cite{Golan}). 
One then defines a semi-ring \emph{ideal} $L$ of $R$ to be simply a ring ideal of $R^\Delta$. 
The \emph{factor semi-ring} $R/L$ is the set of 
equivalence classes in $R$ under the equivalence relation defined by
$$
a \sim b \qquad\text{if and only if}\qquad a-b\in L.
$$
Such a quotient structure will behave as expected. 

Our concern shall be with the path semi-rings $\N \Gamma_+$, for $\Gamma$ an ordinary 1-quiver. 
The elements of $\N\Gamma_+$ are formal linear combinations of \emph{non-empty} paths in $\Gamma$ with
natural co-efficients. Multiplication is given by concatenation of paths.
These rings possess additive cancellation, for $\N \Gamma_+$ embeds into 
the (non-unital) ring $\Z\Gamma_+$. 

\bdf
Let $\Gamma$ be a 1-quiver. An ideal $L\trianglelefteq \N \Gamma_+$ is \textbf{licit} if 
every element $f\colon s\to t$ in $\N \Gamma_+/L$ can be uniquely represented as a sum of edges $s\to t$.
\edf

\bex			\label{X: Licit}
Let $\Gamma$ be the Jordan quiver:
$$
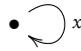
\xymatrix{\bullet\ar@(ur,dr)[]^x} 
$$
An ideal in $\N\Gamma_+$ is licit if 
and only if it is principal and generated by a quadratic polynomial of the form 
$x^2 - a x$, for some $a\in\N$.
\eex

\bth			\label{S: Licit}
Let $\Gamma$ be a 1-quiver with $\Gamma_1$ its set of edges. 
An ideal $L\trianglelefteq \N\Gamma_+$ is licit if and only if the following conditions be met. 
\begin{itemize}
\item The ideal has the form 
$$
L = \Ideal{pq - g_{pq} | p,q\in \Gamma_1}, 
$$
with each $g_{pq}$ a sum of edges (with the same source and target as $pq$). 
\item The obvious compatibility criteria hold: if $g_{pq}=\sum_j s_j$ and $g_{qr}=\sum_k t_k$, then
$$
\sum_j g_{s_j r} = \sum_k g_{p t_k}.
$$
\end{itemize}
\eth

\bpr
If $L$ is licit, define $g_{pq}$ to be the unique expression, modulo $L$, of $pq$ as a sum of edges. Write 
$L'$ for the ideal generated by all elements $pq-g_{pq}$. 
Since each element of $\N\Gamma_+/L'$ can be expressed as a sum of edges, it must be that $L'=L$. 
Moreover, if $g_{pq}=\sum_j s_j$ and $g_{qr}=\sum_k t_k$, then 
$$
\sum_j g_{s_j r} \equiv \sum_j s_j r = g_{pq} r \equiv (pq)r = p(qr) \equiv p g_{qr} = \sum_k p t_k 
\equiv \sum_k g_{p t_k} \ \mod L .
$$
Both sides are sums of edges, and, as each element of $\N\Gamma_+/L$ has a unique such expression, 
it must be that $\sum_j g_{s_j r} = \sum_k g_{p t_k}$.

Conversely, if the two stated conditions are fulfilled, one may endow the monoid $\N\Gamma_1$ (the span of the edges) 
with a multiplication $p\ast q=g_{pq}$. By compatibility, this operation is associative, 
and clearly $\N\Gamma_1\cong \N\Gamma_+/L$.
\epr

\section{2-Quivers}

We now state the definition of a 2-quiver (cf. Subsection 2.2.3 of \cite{Rouquier}).

\bdf				\label{D: 2-quiver}
A \textbf{2-quiver} is a triple $Q=(Q[0],Q[1],Q[2])$, consisting of \textbf{vertices} $Q[0]=\{X,Y,Z,\dots\}$,  
directed \textbf{1-edges} $Q[1]=\{F,G,H,\dots\}$ between vertices, 
and directed \textbf{2-edges}  $Q[2]=\{\alpha,\beta,\gamma,\dots\}$  between 
1-edges having the same source and target. 
At each vertex $X$, a special 1-edge $I_X\colon X\to X$ should be flagged, called the \textbf{stationary 1-edge} 
at $X$.
\edf

All our quivers shall be \emph{finite}, meaning that each $Q[n]$ be a finite set.
Given a 2-quiver $Q=(Q[0],Q[1],Q[2])$, a 1-quiver
$(Q[0],Q[1])$ is obtained from dropping the 2-edges. Moreover, $Q(X,Y)$ is a 1-quiver for any $X,Y\in Q[0]$. 

\bdf				\label{D: Bound}
A \textbf{bound 2-quiver} is a pair $(Q,L)$, where $Q$ is a 2-quiver and $L\trianglelefteq \N Q[0,1]_+$ 
is a licit ideal including the relations 
$$
I_Y F - F \qquad\text{and}\qquad F I_X - F \qquad\text{for any $F\colon X\to Y$.}  
$$
\edf

Abusing terminology, we shall usually drop $L$ from the notation, and refer to just $Q$ as a bound 2-quiver.

\section{Path 1\nicefrac12-Categories}

As an ancillary step towards quiver 2-categories, we shall have reason to consider 2-categories where the horizontal 
composition is ``left undefined''. These gadgets will have 1-cell composition and 
vertical 2-cell composition, and they are required to satisfy all the axioms normally imposed upon a 
2-category, except, of course, those 
pertaining to horizontal composition. 
Let us agree to call such an entity a \textbf{1\nicefrac12-category}.

\bdf 
Let $Q$ be a (bound) 2-quiver. A \textbf{2-path} is a sequence of vertically composable 2-edges 
$\alpha_k\colon F_{k-1}\to F_{k}$, for $k=1,\dots,n$. (Perforce all the 1-cells $F_k$ have the same source and target.)

We denote by $\iota_F\colon F\to F$ the \textbf{empty} or \textbf{stationary} 2-path at $F$.
\edf

\bdf			\label{D: Path Cat}
Let $Q$ be a bound 2-quiver. 
The \textbf{path 1\nicefrac12-category} $\K Q$ is defined in the following way.
\begin{description}
\item[0-cells.] Vertices of $Q$.
\item[1-cells.] Formal finite direct sums of 1-edges of $Q$. 
This includes, for each pair of vertices $X, Y$, a \emph{zero 1-cell} $0\colon X\to Y$.
\item[2-cells.] Matrices of linear combinations of 2-paths in $Q$.
\item[1-cell composition.] Concatenation of 1-edges, modulo the binding relations of $Q$. 
By the definition of licit ideal, 
this operation is well defined and associative. 
\item[Vertical 2-cell composition.] Concatenation of 2-paths, combined with matrix multiplication. 
\end{description}
\edf

The path 1\nicefrac12-category has a natural gradation by path lengths, and one may write
$$ 
\K Q = \bigoplus_{n=0}^\infty \K Q_n,
$$
where $Q_n$ stands for the set of paths of length $n$.
 
\bdf
The \textbf{edge ideal} 
$$
\K Q_+ = \bigoplus_{n=1}^\infty \K Q_n
$$ 
is the vertical ideal of $\K Q$ generated by the 2-edges of $Q$.  
\edf 

Since the finitary property places no restrictions on horizontal composition, 
the definition carries over unaltered, so that one may freely speak of finitary 1\nicefrac12-categories. 

\bth	
Let $Q$ be a bound 2-quiver.
The bound path 1\nicefrac12-category $\K Q$ has the following properties: 
\bnum
\item It has finitely many 0-cells. 
\item It has fully additive 1-categories $\K Q(X,Y)$. 
\item It has finitely many indecomposable 1-cells, namely the 1-edges of $Q$.
\item The identity 1-cells are indecomposable. 
\item The vector spaces $\K Q(F,G)$ are finite-dimensional if and only if 
there are no oriented cycles on the 2-level in $Q$.
\enum
Hence $\K Q$ is finitary if and only if there are no oriented cycles on the 2-level in $Q$.
\hfill
\eth

\bpr
Claim 1 is obvious. Claims 3 and 4 follow from the fact that, by the definition of licit ideal, the indecomposable 1-cells are precisely the 1-edges. Claim 5 follows from the obvious fact that $Q$ contains finitely many 2-paths if and only if it has no cycles on the 2-level. 

To prove claim 2, we note that $\K Q(X,Y)$ is equivalent to the fully additive category of projective right modules over the algebra 
$\bigoplus \K Q(F,G)$, the summation extending over all 1-edges $F,G\colon X\to Y$.
This completes the proof.
\epr

\bex 			\label{X: Quiver}
Consider the bound 2-quiver given pictorially:
\vspace{-3ex}
$$ 
\xymatrix{ 
X \ruppertwocell<10>^{I_X}{\alpha} \ar[r]|{F} \rlowertwocell<-10>_{I_X}{\beta} & X
}
\qquad F^2=I_X\oplus F.
\vspace{-3ex}
$$
The path 2-category comprises a single 0-cell $X$ and 
\begin{itemize}
\item 1-cells, e.g.: $0$, $I_X$, $F$, $I_X\oplus F$, $I_X\oplus F\oplus F$, \&c
\item 2-cells, e.g.: 
\begin{gather*}
\iota_{I_X}\colon I_X \to I_X, \qquad \iota_F\colon F\to F, \qquad \alpha\colon I_X\to F, \\
 \alpha\circ\beta\colon F\to F, \qquad \iota_{I_X}+2\beta\circ\alpha\colon I_X \to I_X, \\ 
\begin{pmatrix}
\iota_{I_X}+2\beta\circ\alpha & 0 \\ 
-\alpha & \alpha\circ\beta 
\end{pmatrix} \colon I_X\oplus F \to I_X\oplus F, \qquad \text{\&c.}
\end{gather*}
\end{itemize}
It is not finitary because $\alpha$ and $\beta$ form a closed cycle.
\eex

\section{Admissible Ideals}

Let $C$ be a 2-category, or even 1\nicefrac12-category, and let $J$ be a vertical ideal. 
 For a positive integer $n$, 
we define the power $J^n$ as the vertical ideal generated by vertical $n$'th powers of 2-cells in $J$.

\bdf 			\label{D: Admissible}
Let $Q$ be a bound 2-quiver. A vertical ideal $J\trianglelefteq \K Q$ is \textbf{admissible} if, 
for some $n$,
$$
\K Q_+^n\subseteq J\subseteq \K Q_+^2.
$$
The quotient $\K Q/J$ will be called a \textbf{quiver 1\nicefrac12-category}.
\edf

\bth
Admissible ideals are finitely generated as bimodules over $\K Q$. 
\eth

\bpr
For each pair $X$, $Y$ of 0-cells (finitely many), there is an exact sequence of bimodules:
$$ 
\xymatrixcolsep{2pc}
\xymatrix{0\ar[r] &  \K Q_+^n(X,Y) \ar[r] & J(X,Y) \ar[r] & J(X,Y)/\K Q_+^n(X,Y) \ar[r] & 0 }
$$
The bimodule on the left is generated by the finitely many 2-paths of length $n$ between indecomposable 1-cells.  
The bimodule on the right embeds in 
$\K Q/\K Q_+^n$, which is finitely generated since $\K Q$ 
only contains finitely many paths, between indecomposable 1-cells, of length not exceeding $n$. 
Consequently, the middle term is finitely generated. 
\epr

\section{Vertical Drop}

Let $Q$ be a bound 2-quiver supplied with an admissible ideal $J$, and consider the 
1\nicefrac12-category $\K Q/J$. Introduce the notation:
\begin{align*}
\K Q/J[X,-] & \quad \text{for the set of 2-cells joining 1-cells with source $X$;} \\
\K Q/J[-,Y] & \quad \text{for the set of 2-cells joining 1-cells with target $Y$.}
\end{align*}

\bdf 			\label{D: Drop}
Suppose given, for each 1-cell $F\colon X\to Y$ in $\K Q/J$, two linear transformations 
$$
F_\ast\colon \K Q/J[-,X]\to \K Q/J[-,Y] \qquad\text{and}\qquad F^\ast\colon \K Q/J[Y,-]\to \K Q/J[X,-], 
$$
both operators $(-)_\ast$ and $(-)^\ast$ being additive and subject to the following axioms: 
$$ \setlength\arraycolsep{2pt}
\begin{array}{rcrclrcl} 
\text{\textsc{i}.} &\quad& (FG)_\ast&=&F_\ast G_\ast & (FG)^\ast&=& G^\ast F^\ast \\
\text{\textsc{ii}.} && (I_X)_\ast &=& 1_{\K Q/J[-,X]} & (I_X)^\ast &=& 1_{\K Q/J[X,-]} \\
\text{\textsc{iii}.} && F_\ast(\alpha\circ\beta) &=& F_\ast(\alpha)\circ F_\ast(\beta) & F^\ast(\alpha\circ\beta) &=& F^\ast(\alpha)\circ F^\ast(\beta) \\
\text{\textsc{iv}.} && F_\ast(\iota_K) &=& \iota_{FK} & F^\ast(\iota_K) &=& \iota_{KF} \\
\text{\textsc{v}.} &&&& \multicolumn{2}{c}{F_\ast G^\ast =  G^\ast F_\ast} \\
\text{\textsc{vi}.} &&&& 
\multicolumn{2}{c}{\tilde{F}_\ast(\beta)\circ G^\ast(\alpha)  =  \tilde{G}^\ast(\alpha)\circ F_\ast(\beta)}
\end{array}
$$
For Axiom \textsc{vi}, we assume ourselves placed in the situation described pictorially by:
\beq			\label{E: Horizontal}
\xymatrix{
Z \ruppertwocell^{G}{\beta} \ar[r]_{\tilde{G}} 
& Y \ruppertwocell^{F}{\alpha} \ar[r]_{\tilde{F}} 
& X
}
\eeq
Then these maps $(-)_\ast$ and $(-)^\ast$ will be referred to as \textbf{vertical drop relations}.
\edf 

\bth			\label{S: Horizontal}
Let $Q$ be a bound 2-quiver with an admissible ideal $J$, and 
suppose the 1\nicefrac12-category $\K Q/J$ has been endued with vertical drop relations as above. 
Defining an horizontal composition 
$$
\alpha\beta =  \tilde{G}^\ast(\alpha)\circ F_\ast(\beta) = \tilde{F}_\ast(\beta)\circ G^\ast(\alpha)
$$
in the situation represented by \eref{E: Horizontal}, will then make $\K Q/J$ into a 2-category.  
\eth

\bpr
The horizontal composition thus defined on $\K Q/J$ is associative, for, assuming the 
situation pictorially represented by: 
$$
\xymatrix{
W \ruppertwocell^{H}{\gamma} \ar[r]_{\tilde{H}} 
& Z \ruppertwocell^{G}{\beta} \ar[r]_{\tilde{G}} 
& Y \ruppertwocell^{F}{\alpha} \ar[r]_{\tilde{F}} 
& X
}
$$
one has 
\begin{align*}
(\alpha\beta)\gamma &= \tilde{H}^\ast( \tilde{G}^\ast(\alpha)\circ F_\ast(\beta) ) \circ (FG)_\ast(\gamma) \\
&= \tilde{H}^\ast \tilde{G}^\ast(\alpha) \circ \tilde{H}^\ast F_\ast(\beta)  \circ F_\ast G_\ast(\gamma) \\
&= \tilde{H}^\ast\tilde{G}^\ast(\alpha) \circ F_\ast\tilde{H}^\ast(\beta)\circ F_\ast G_\ast(\gamma) \\
&= (\tilde{G}\tilde{H})^\ast(\alpha) \circ F_\ast(\tilde{H}^\ast(\beta)\circ G_\ast(\gamma)) 
= \alpha(\beta\gamma) ,
\end{align*}
using Axioms \textsc{i}, \textsc{iii} and \textsc{v}.

Let now $\alpha\colon F\to \tilde F$ for $F,\tilde F\colon Y\to X$. Then, using Axioms \textsc{ii} and \textsc{iv}, we find 
\begin{align*}
\alpha\iota_{I_Y} &= I_Y^\ast(\alpha)\circ F_\ast(\iota_{I_Y}) = \alpha\circ \iota_{FI_Y} = \alpha  \\
\iota_{I_X}\alpha &= \tilde{F}^\ast(\iota_{I_X})\circ (I_X)_\ast(\alpha) = \iota_{I_X\tilde{F}}\circ\alpha 
= \alpha, 
\end{align*}
 establishing that horizontal composition is unital. Moreover, 
$$
\iota_F\iota_G = G^\ast(\iota_F)\circ F_\ast(\iota_G) = \iota_{FG}\circ\iota_{FG} = \iota_{FG}
$$
by Axiom \textsc{iv}.

Finally, let us verify the Interchange Law. Consider the situation represented by the following diagram:
$$ 
\xymatrix{
\ast \ruppertwocell^{\tilde{H}}{\tilde\beta} \ar[r]|(0.75){\tilde G} \rlowertwocell_{\tilde F}{\tilde\alpha}
& \ast \ruppertwocell^H{\beta} \rlowertwocell_F{\alpha} \ar[r]|(0.75)G & \ast
}
$$ 
Using Axioms \textsc{iii} and \textsc{vi}, we find 
\begin{align*}
\alpha \tilde\alpha \circ \beta \tilde\beta &= \tilde F^\ast(\alpha) \circ G_\ast(\tilde\alpha) 
\circ \tilde G^\ast(\beta) \circ H_\ast(\tilde\beta) \\
&= \tilde F^\ast(\alpha)\circ\tilde F^\ast(\beta) \circ H_\ast(\tilde\alpha) \circ H_\ast(\tilde\beta) \\
&= \tilde F^\ast(\alpha\circ\beta) \circ H_\ast(\tilde\alpha \circ \tilde\beta)
= (\alpha\circ\beta)(\tilde\alpha \circ \tilde\beta) . \qedhere
\end{align*} 
\epr

\section{Quiver 2-Categories} 	\label{A: Quiver 2-Categories}

\bdf			
A 2-category obtained in the way described in Theorem \ref{S: Horizontal} will be called a \textbf{quiver 2\hyp category}.
\edf

\bth
The Jacobson radical of a quiver 2-category $\K Q/J$ is  
$\K Q_+/J$.
\eth

\bpr
Since $\K Q_+^n\subseteq J$, the vertical ideal $\K Q_+/J$ will be a (vertically) nilpotent 
ideal of $\K Q/J$. Moreover, the quotient
$$ 
\frac{ \K Q/J}{ \K Q_+/J } \cong \K Q/\K Q_+ 
$$
is discrete. Therefore $\Rad \K Q/J = \K Q_+/J$. 
\epr

\bth			\label{S: Finitary}
Quiver 2-categories are finitary, the indecomposable 1-cells being precisely the 1-edges of the underlying 2-quiver.
\eth

\bpr 
Obviously, $\K Q/J$ contains finitely many objects, and hence condition 1 in Definition~\ref{D: Finitary} is satisfied.

Next, we turn to conditions 3 and 4. We show that the indecomposable 1-cells of $\K Q/J$ coincide with those of $\K Q$. 
Since $J$ contains no identity 2-cells $\iota_F$, all 1-cells of $\K Q$ survive when passing to the quotient $\K Q/J$. Clearly, 
a decomposable 1-cell of $\K Q$, say $F=G_1\oplus G_2$, cannot be made indecomposable unless one of its 
constituents $G_k$ be killed. Conversely, 
since the defining relations $\pi_i\lambda_i=\iota_{F_i}$ and $\sum_i\lambda_i\pi_i=\iota_F$ 
of a direct sum system $F=F_1\oplus F_2$ (projections $\pi_i$; injections $\lambda_i$) 
involve identity 2-cells, which are of degree $0$, 
they cannot belong to $J$, so that an indecomposable 1-cell of $\K Q$ will not decompose in $\K Q/J$. 
Similarly, no new isomorphisms can be created between two indecomposable 1-cells of $\K Q$. 

Condition 2 follows since $\K Q/J(X,Y)$ is equivalent to the fully additive category of projective right modules over the algebra 
$\bigoplus \K Q/J(F,G)$, the summation extending over all 1-edges $F,G\colon X\to Y$.
%
%

Finally, condition 5 is satisfied as there are only finitely many paths of 
length not exceeding $n$, and all longer paths vanish in the quotient $\K Q/J$.
\epr

\section{The Gabriel 2-Quiver}		\label{A: Gabriel 2-Quiver}

We now purpose to study a given finitary 2-category $C$ by means of a certain associated 2-quiver.  

\bdf 
The \textbf{(bound) Gabriel 2-quiver} $Q_C$ of a finitary 2-category $C$ is defined as follows. 
\begin{description}
\item[Vertices.] A fixed collection of representatives of equivalence classes of 0-cells of $C$. 
\item[1-edges.] A fixed collection of representatives of isomorphism classes of the indecomposable 1-cells of $C$.
\item[2-edges.] Given two indecomposable 1-cells $F,G\colon X\to Y$ of $C$, 
the number of  2-edges from $F$ to $G$ is declared to be
 $$ \dim \frac{\Rad C}{\Rad^2 C}(F,G).  $$ 
\item[Binding ideal.] The Gabriel quiver is bound by the ideal $\Ker \Xi$, where 
$$\Xi\colon \N Q_C[0,1]\to C[0,1]$$ is the unique 2-functor 
extending the identity functor on our fixed collection of representatives of the 1-cells of $C$. (This ideal is licit as a consequence of the Krull--Schmidt property of $C$.)
\end{description}
\edf 

The Gabriel 2-quiver of a finitary 2-category is clearly finite. 
It is also \textbf{reduced} in the sense that no two 0-cells are equivalent (modulo the binding ideal).

\bth			\label{S: Return}
Bi-equivalent 2-categories have the same Gabriel 2-quiver (up to isomorphism), and
taking the Gabriel quiver of a reduced quiver 2-category returns the original 2-quiver it was based upon. 
\eth

\bpr
The first assertion is obvious.
For the second assertion, the proceedings on levels 0 and 1 are clear. On the 2-level, one has 
$$
\frac{\Rad \K Q/J}{\Rad^2 \K Q/J} = \frac{\K Q_+/J}{\K Q_+^2/J} \cong \K Q_+/\K Q_+^2,
$$
and the theorem follows. 
\epr

\bth			\label{S: Every}
Every finitary 2-category is bi-equivalent with a quiver 2-category for a unique underlying reduced 2-quiver. 
\eth

\bpr
Let $C$ be finitary. Choose a collection $X_1,\dots,X_k$ of representatives of equivalence classes of 0-cells, and then 
a collection $F_1,\dots,F_m$ of representatives of isomorphism classes of the indecomposable 1-cells of the categories 
$C(X_i,X_j)$. Letting $C'$ denote the full sub-2-category of $C$ on the 1-cells $F_1,\dots,F_m$, 
the 2-edges of the Gabriel quiver $Q$ of $C$ will be in one-to-one correspondence with a basis 
$\overline\alpha_1,\dots,\overline\alpha_n$ for the space $\Rad C'/\Rad^2 C'$. 
Because $\K Q$ is free, we may define an obvious 2-functor $\Psi\colon \K Q\to C'$.

Every element of $\Rad C'$ can be expressed as a vertical polynomial in the 2-cells $\alpha_1,\dots,\alpha_n$ 
(assuming the base field $\K$ is algebraically closed; see Section \ref{A: Prel}). 
This establishes the fullness of $\Psi$ on the 2-level. 
By construction, $\Psi$ is surjective on levels 0 and 1. Hence $\Psi$ provides a bi-equivalence   
$$ 
C\cong C'= \Im \Psi \cong \K Q /\Ker \Psi,
$$
where the kernel $\Ker \Psi$ is the set of 2-cells annihilated by $\Psi$.

Next, it is clear that $\Psi(\K Q_+)\subseteq\Rad C$. 
Since $\Rad C$ is (vertically) nilpotent, there exists an $n$ such that 
$$  \Psi(\K Q_+^n)\subseteq \Rad^n C = 0, $$
implying $ \K Q_+^n\subseteq \Ker\Psi$.  

Suppose $\eta+\zeta\in \Ker\Psi$, where $\eta$ is a linear combination of 
the 2-edges $\alpha_1,\dots,\alpha_n$ and $\zeta\in\K Q_+^2$. 
Then $\Psi(\zeta)\in\Rad^2 C$, so
$$ 
0 = \Psi(\eta+\zeta) \equiv \Psi(\eta) \ \mod \Rad^2 C.
$$ 
Since $\alpha_1,\dots,\alpha_n$ are linearly independent modulo $\Rad^2 C$, 
it must be that $\eta=0$, and so $\eta+\zeta=\zeta\in\K Q_+^2$. 
Consequently, $\Ker\Psi\subseteq \K Q_+^2$, and we have shown that $\Ker\Psi$ is admissible.

There remains to check the vertical drop relations. Given a 1-cell $F$, define (when applicable) 
$$
F_\ast(\xi) = \iota_F\xi \qquad\text{and}\qquad F^\ast(\xi) = \xi\iota_F. 
$$
The first five axioms for vertical drop are immediate consequences of the properties of horizontal composition. 
Axiom \textsc{vi} follows from the calculation 
\begin{align*}
\tilde{F}_\ast(\beta)\circ G^\ast(\alpha)  &= \iota_{\tilde F}\beta\circ\alpha\iota_G = 
( \iota_{\tilde F}\circ\alpha)( \beta\circ\iota_G) = \alpha\beta \\
&= (\alpha\circ\iota_F)(\iota_{\tilde G} \circ \beta) 
= \alpha\iota_{\tilde G} \circ \iota_F\beta =\tilde{G}^\ast(\alpha)\circ F_\ast(\beta),
\end{align*} 
which also shows that these vertical drop relations will indeed induce the original horizontal composition on 
$C\cong\K Q/J$.

Uniqueness of the 2-quiver follows directly from Theorem \ref{S: Return}. 
\epr

\section{Examples}

We close by presenting a sketch-book of examples.

\bex
The following method for constructing a 2-category out of an ordered monoid was devised by Kudryavtseva \& 
Mazorchuk (unpublished). 

A partial ordering $\leq$ on a monoid is said to be
\emph{compatible} with the multiplication if $x\leq y$ implies $xz\leq yz$ and $zx\leq zy$ 
for all $x,y,z\in M$. Let $M$ be a finite monoid with identity element $1$ and a
 compatible ordering $\leq$. 

Define the finitary 2-category $C_M$ as follows. The single 0-cell is $\ast$. The 
indecomposable 1-cells are the elements of $M$. 
Composition of 1-cells is monoid multiplication. Next, for each relation $x\leq y$, introduce a formal symbol $e_{x,y}$ and define
$$
C_M(x,y) = 
\begin{cases}
\K e_{x,y} & \text{if $x\leq y$} \\
0 & \text{else.}
\end{cases}
$$
Vertical composition is given by the rule $e_{y,z}\circ e_{x,y}=e_{x,z}$, and horizontal composition by 
$e_{u,v}e_{x,y}=e_{ux,vy}$ (both extended by bilinearity). 

The 2-category $C_M$ may be encoded as follows. Its Gabriel 2-quiver $Q$ has a single vertex $\ast$ 
and 1-edges $F_x$, corresponding to the elements $x\in M$. The radical of $C_M$ is generated by all 
$e_{xy}$, where $y$ covers $x$ in the ordering $\leq$,
which means there should be a 2-edge $\epsilon_{x,y}$ in $Q$ for each $y$ \emph{covering} $x$. 
($y$ is said to cover $x$ if $x<y$ and there is no $x<z<y$.)
In short, $Q$ is the Hasse diagram of the partially ordered set $M$. 
It is evidently bound by relations $F_xF_y=F_{xy}$, for $x,y\in M$. 

The admissible ideal is generated by all differences of converging paths: elements of the form 
$$
\epsilon_{y_m,z} \circ\cdots\circ \epsilon_{y_1,y_2} \circ \epsilon_{x,y_1} - 
\epsilon_{\tilde y_n,z} \circ\cdots\circ \epsilon_{\tilde y_1,\tilde y_2} \circ \epsilon_{x,\tilde y_1}.
$$ 
 Finally, the vertical drop relations will be found from an inspection of $M$. 
For example, to determine 
$(F_z)_\ast(\epsilon_{x,y})$,
one has to express $e_{z,z}e_{x,y}=e_{zx,zy}$ as a path in the 2-edges $\epsilon$. 
\eex

\bex				
As a special case of the preceding example, we may consider Grensing \& Mazorchuk's 
categorification of the Catalan monoid 
in \cite{Catalan}. According to their Theorem 1,  
the $n$'th Catalan monoid $C_n$ of all order\hyp preserving and order-decreasing full transformations on $\{1,2,\dots,n\}$ 
(it can also be realised as a Hecke--Kiselman monoid of type $A_{n-1}$) 
is categorified by a 2-category $D_n$ having a single
0-cell $\ast$ and as 1-cells the elements of $C_n$. 
The 2-cells are
$$  
D_n(p,q)=
\begin{cases}
\K\pi_{pq} & \text{if $p\leq q$} \\
0 & \text{else.}
\end{cases}
$$ 
Here $\leq$ is the natural ordering on $C_n$, and $\pi_{pq}$, whenever $p\leq q$, denotes 
a certain ``natural projection''. (See Proposition 9 of \cite{Catalan} for details.)

Since the ordering $\leq$ is compatible with the multiplication in $C_n$,  
the Gabriel quiver may be constructed as in the preceding example.  
\eex

\bex			\label{X: S_A}
A method for constructing 2-categories using projective functors on module categories of
finite-dimensional algebras was proposed by Mazorchuk \& Miemietz in Subsection 7.3 of \cite{Cell}. 
We describe here a slightly modified variant, which may be contrasted with the original version 
to be found below in Example \ref{X: C_A}.

Let $A$ be a finite-dimensional, connected, basic algebra (associative and unital as always), arising from 
a connected Gabriel 1-quiver $\Gamma_A$ and an admissible ideal $K$ with $\K\Gamma_A/K\cong A$. 
We define what might be called a \emph{semi-2-category} $S_A$, there being no identity 1-cells, as follows. 
There should be a single 0-cell~$\ast$, identified with (a small category equivalent to) the module category $\AMod$, 
the indecomposable 1-cells being the so-called \emph{projective functors}. These are functors isomorphic to 
tensoring with projective $A$--$A$-bimodules. 2-cells are natural transformations of functors. 
As there are only finitely many indecomposable $A$--$A$-bimodules up to isomorphism, $S_A$ is finitary. 
 
Alternatively, one may construct a bi-equivalent semi-2-category 
(which we, by abuse of notation, also may denote by $S_A$), as follows.
It has a single 0-cell $\ast$ and its 1-cells are direct sums of bimodules isomorphic with projective $A$--$A$-bimodules. 
Such a module is of the form $Ae\otimes f A$, with $e$ and $f$ being primitive, orthogonal idempotents of $A$. 
Composition of 1-cells is 
the tensor product of bimodules, obeying the rule
$$
(Ae \otimes fA) \otimes_A (Ag \otimes hA) = Ae \otimes fAg \otimes hA = 
(Ae \otimes hA)^{\oplus \dim fAg}.
$$ 
2-cells are module homomorphisms, the horizontal composition being the tensor product of bimodules and 
vertical composition the usual composition of homomorphisms.
An homomorphism 
$$ 
\phi\colon Ae\otimes fA \to Ag\otimes hA
$$
is specified by a pair of (linear combinations of) paths in $\Gamma_A$:  
$$
\phi=(p\colon g\to e,\ q\colon f\to h),\quad 
e\otimes f \mapsto p\otimes q. 
$$

The Gabriel 2-quiver $Q$ of $S_A$ will contain a single vertex $\ast$ 
and 1-edges $F_{ef}$ for each pair $(e,f)$ of vertices in $\Gamma_A$,  bound by relations 
$$ 
F_{ef} F_{gh} = F_{eh}^{\oplus \dim fAg}.
$$ 
The radical of $S_A$ is generated by pairs $(x,1)$ and $(1,y)$, where 
$x$ and $y$ are edges of $\Gamma_A$. Therefore, the 2-edges of $Q$ are
$$
\alpha_{g,x} \colon F_{ge} \to F_{gf} \qquad\text{and}\qquad \alpha_{x,g}\colon F_{fg} \to F_{eg},
$$
for all edges $x\colon e\to f$ and vertices $g$ in $\Gamma_A$. In other words, 
$Q[1,2]\cong \Gamma_A\times\Gamma_A^\circ$, where $\Gamma_A^\circ$ denotes the dual quiver. 
The principal reason for this phenomenon is that the quiver just mentioned encodes homomorphisms 
between indecomposable projective $A$--$A$-bimodules, that is, $A\otimes A^\circ$-modules. 

The admissible relations on $Q$ take on either of two shapes. 
\bnum
\item Relations inherited from $K$:
$$ 
p(\alpha_{g,x_1},\dots,\alpha_{g,x_n}) = p(\alpha_{x_1,h},\dots,\alpha_{x_n,h}) = 0 
$$ 
for all vertices $g,h$ and each polynomial relation $p(x_1,\dots,x_n)\in K$.
\item Commutation relations 
$$ \alpha_{g,x}\circ \alpha_{y,e} = \alpha_{y,f}\circ \alpha_{h,x}, $$  
for any edges $x\colon e\to f$ and $y\colon g\to h$.
\enum

As a concrete illustration, we take the path algebra $A$ on the quiver: 
$$
\xymatrix{e \upar[r]^x & f } 
$$
The Gabriel 2-quiver of $S_A$ is: 
$$
\xymatrix{
F_{fe} \ar[r]^{\alpha_{f,x}} \ar[d]_{\alpha_{x,e}}  
& F_{ff} \ar[d]^{\alpha_{x,f}}  \\
F_{ee} \ar[r]_{\alpha_{e,x}} & F_{ef}
} 
$$
It is bound by equations 
$$
F_{ab}F_{cd} = 
\begin{cases}
F_{ad} & \text{if $(b,c)=(e,f)$} \\
0 & \text{else;}
\end{cases}
$$
and a single admissible relation 
$$
\alpha_{e,x}\circ \alpha_{x,e} = \alpha_{x,f}\circ \alpha_{f,x}. \qedhere
$$
\eex

\bex				\label{X: C_A}
The original version $C_A$ of the category defined in the previous example also incorporated an
identity 1-cell on $\ast$, making $C_A$ a true 2-category.
This is the functor $I=A\otimes_A -$, corresponding to the 
bimodule $A$. Its inclusion will unfortunately destroy the 2-quiver just described, for which the erratic 
homomorphisms into and out of $A$ are to blame. 

Again we illustrate with the path algebra $A$ on the quiver: 
$$
\xymatrix{e \upar[r]^x & f } 
$$
A sketch of the category $C_A$, displaying a basis for the spaces of homomorphisms, looks as follows: 
$$
\xymatrix{
Af\otimes eA \ar[rr]^{f\otimes x} \ar[dd]_{x\otimes e} \ar[dr]^{x} 
&& Af\otimes fA \ar[dd]^{x\otimes f} \ar[dl]^{f} \\
& A \ar[dr]_{x\otimes f + e\otimes x} \\
Ae\otimes eA \ar[rr]_{e\otimes x} \ar[ur]_{e}  && Ae\otimes fA
} 
$$
Each bimodule is monogenic, and the symbol on each arrow should be interpreted as the image of 
the obvious generator.

The resulting Gabriel 2-quiver is: 
$$
\xymatrix{
F_{fe} \ar[rr]^{\gamma_2} \ar[dd]_{\gamma_1}  
&& F_{ff}  \ar[dl]^{\beta_2} \\
& I \ar[dr]_{\alpha} \\
F_{ee} \ar[ur]_{\beta_1}  && F_{ef}
} 
$$
It is bound by relations  
$$
F_{ab}F_{cd} = 
\begin{cases}
F_{ad} & \text{if $(b,c)=(e,f)$} \\
0 & \text{else;}
\end{cases}
$$
and an admissible ideal generated by 
$$
\beta_1\circ\gamma_1 = \beta_2\circ \gamma_2.
$$
One may contrast this with the quiver obtained for $S_A$ in Example \ref{X: S_A}.
\eex

\bex					\label{X: Duflo}
Our last example constructs the first known example of a fiat 2-category possessing a Duflo 
involution which is not self-adjoint. Fiat 2-categories, introduced by Mazorchuk \& Miemietz, are
finitary 2-categories with a weak involution $\star$ (anti\hyp equivalence of order 2) and adjunction morphisms
(between $F$ and $F^{\star}$); see \cite{Cell} for details. 
The combinatorial structure of finitary categories is encoded in the so-called \emph{left}, \emph{right} and
\emph{two-sided cells}. For example, two indecomposable 1-cells $F$ and $G$ belong to the same left
cell provided that there be 1-cells $H$ and $K$ such that $F$ is isomorphic to a summand of $HG$
and $G$ is isomorphic to a summand of $KF$. In a similar manner, one defines right and two-sided cells.
  
In \cite{Cell}, it was shown that each left cell of a finitary 2-category
contains a special 1-cell called its \emph{Duflo involution}. This 1-cell
plays a crucial rôle in the construction of cell 2-representations of fiat 2-categories.
The examples considered in \cite{Cell} share the property that the Duflo involutions of left cells
are self-dual with respect to $\star$. 
We here present a fiat category with a Duflo involution that is \emph{not} self-dual.

Consider the path algebra $V$ based upon the cyclic Kronecker quiver:  
$$
\xymatrix{
e \upar[r]^x & f \upar[l]^y
} \qquad xy=yx=0.
$$
We construct the finitary 2-category $C_V$ as described in Example \ref{X: C_A}.
The indecomposable 1-cells are $V$ itself, along with the projective bimodules: 
$$
P_1 = Ve\otimes eV, \qquad P_2 = Vf\otimes fV, \qquad Q_1 = Vf\otimes eV, \qquad Q_2 = Ve\otimes fV.  
$$
The composition table is given in Table \ref{T: 1-comp}.
\begin{table}
$$
\begin{array}{c|cccc}
\circ & P_1 & P_2 & Q_1 & Q_2 \\
\hline
P_1 & P_1 & Q_2 & P_1 & Q_2 \\
P_2 & Q_1 & P_2 & Q_1 & P_2 \\
Q_1 & Q_1 & P_2 & Q_1 & P_2 \\
Q_2 & P_1 & Q_2 & P_1 & Q_2
\end{array}
$$
\caption{1-Cell Composition in $C_V$.\label{T: 1-comp}}
\end{table}
The homomorphisms of these bimodules are, apart from identity transformations, 
(vertically) generated by those indicated in the diagram:
$$ 
\xymatrix{ 
&& \downar[dll]_{\alpha_{11}}  P_1 \upar[d]^{\gamma_1} \upar[drr]^{\alpha_{12}} \\
Q_1 \downar[urr]_(0.75){\beta_{11}} \upar[drr]^(0.75){\beta_{12}} \upar[rr]^{\delta_1} 
&& \upar[ll]^{\eta_1} V \downar[d]_{\zeta_2} \upar[u]^{\zeta_1} \downar[rr]_{\eta_2} 
&& \downar[dll]_(0.75){\beta_{22}} Q_2 \downar[ll]_{\delta_2} \upar[ull]^(0.75){\beta_{21}} \\
&& \upar[ull]^{\alpha_{21}}  P_2 \downar[u]_{\gamma_2} \downar[urr]_{\alpha_{22}} 
} 
$$
Rules for vertical and horizontal composition are provided in Tables \ref{T: V-comp} and~\ref{T: H-comp}, 
respectively.
\begin{table} 
\begin{align*}
\alpha_{jk} \circ \beta_{ij} = \alpha_{j'k} \circ \beta_{ij'} &= \eta_k \circ \delta_i
& \beta_{jk} \circ \alpha_{ij} = \beta_{j'k} \circ \alpha_{ij'} &= \zeta_k \circ \gamma_i \\
\alpha_{ji} \circ \beta_{ij} &= 0  & \beta_{ji} \circ \alpha_{ij} &= 0 \\
\alpha_{ij} \circ \zeta_i &= 0 	& 	\gamma_j \circ \beta_{ij} &= \delta_i \\
\delta_j \circ \alpha_{ij} &= 0 	&  \beta_{ij} \circ \eta_i &= \zeta_j \\
\gamma_i \circ \zeta_i = \zeta_i \circ \gamma_i &= 0 & \delta_i\circ \eta_i = \eta_i \circ \delta_i &= 0 \\
\zeta_j \circ \delta_i &= 0 & \eta_j \circ \gamma_i &= \alpha_{ij}
\end{align*}
\caption{Vertical Composition in $C_V$.\label{T: V-comp}}
\end{table}
(We have chosen to specify horizontal compositions with $\iota_{P_1}$ and $\iota_{Q_1}$ only. 
To obtain products with $\iota_{P_2}$ and $\iota_{Q_2}$, one takes advantage of the inherent symmetry 
between indices $1$ and $2$.) 
\begin{table}
$$
\begin{array}{c|ccccc}
& (\iota_{P_1})_\ast & (\iota_{P_1})^\ast & (\iota_{Q_1})_\ast & (\iota_{Q_1})^\ast \\
\hline
\alpha_{11} & \iota_{P_1} & \alpha_{11} & \iota_{Q_1} & \alpha_{11} \\
\alpha_{12} & \alpha_{12} & \iota_{P_1} & \beta_{12} & 0_{P_1} \\
\alpha_{21} & \beta_{21} & 0_{Q_1} & \alpha_{21} & \iota_{Q_1} \\
\alpha_{22} & 0_{Q_2} & \beta_{11} & 0_{P_2} & \beta_{11} \\
\beta_{11} & 0_{P_1} & \beta_{11} & 0_{Q_1} & \beta_{11} \\
\beta_{12} & \alpha_{12} & \iota_{Q_1} & \beta_{12} & 0_{Q_1} \\
\beta_{21} & \beta_{21} & 0_{P_1} & \alpha_{21} & \iota_{P_1} \\
\beta_{22} & \iota_{Q_2} & \alpha_{11} & \iota_{P_2} & \alpha_{11} \\
\gamma_1 & \iota_{P_1} & \iota_{P_1} & \iota_{Q_1} & \alpha_{11} \\
\gamma_2 & \beta_{21} & \beta_{11} & \alpha_{21} & \iota_{Q_1} \\
\delta_1 & 0_{P_1} & \beta_{11} & 0_{Q_1} & 0_{Q_1} \\
\delta_2 & \beta_{21} & 0_{P_1} & \alpha_{21} & \alpha_{11} \\
\zeta_1 & 0_{P_1} & 0_{P_1} & 0_{Q_1} & \beta_{11} \\
\zeta_2 & \alpha_{12} & \alpha_{11} & \beta_{12} & 0_{Q_1} & \\
\eta_1 & \iota_{P_1} & \alpha_{11} & \iota_{Q_1} & \iota_{Q_1} \\
\eta_2 & \alpha_{12} & \iota_{P_1} & \beta_{12} & \beta_{11}
\end{array}
$$
\caption{Horizontal Composition in $C_V$.\label{T: H-comp}}
\end{table}

Now consider the bound 2-quiver (with a single 0-cell):
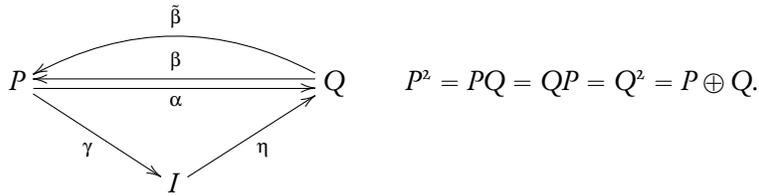
$$
\xymatrix{
P \downar[rr]_{\alpha} \ar[dr]_{\gamma}
&& Q \downar[ll]_{\beta} \cruvear[ll]_{\tbeta} \\
& I  \ar[ur]_{\eta}
}
\qquad P^2=PQ=QP=Q^2=P\oplus Q. 
$$
Construct a quiver 2-category $M$ by factoring out the admissible relations: 
\begin{align*}
\beta\circ\alpha &= 0 & \alpha\circ\beta &= 0 \\ 
\gamma\circ \beta &= \gamma\circ\tbeta & \beta\circ\eta &= \tbeta\circ\eta \\ 
\beta\circ\eta\circ\gamma &= \tbeta\circ\alpha   & \eta\circ\gamma\circ\beta  &= \alpha\circ\tbeta \\
&& \gamma\circ\beta\circ\eta &= 0 ,
\end{align*}
along with the vertical drop relations in Table \ref{T: M Drop}.  
\begin{table}
$$
\begin{array}{c|cccc}
		& P_\ast & P^\ast & Q_\ast & Q^\ast \\
		\hline
\alpha & \iota_P & \alpha+\beta & \iota_Q & \alpha+\beta \\ 
\beta  & \iota_Q & \alpha+\beta & \iota_P & \alpha+\beta \\ 
\tbeta & -\alpha+\tbeta+\eta\circ\gamma & \iota_Q & -\alpha+\tbeta+\eta\circ\gamma & \iota_P \\
\gamma & \iota_P + \tbeta & \iota_P + \beta & \iota_Q -\alpha + \eta\circ\gamma & \iota_Q + \alpha \\ 
\eta & \iota_P  -\alpha+\eta\circ\gamma & \iota_P + \alpha & \iota_Q + \tbeta & \iota_Q + \beta
\end{array}
$$
\caption{Vertical Drop relations for $M$.\label{T: M Drop}}
\end{table}
One should think of $M$ as modelled on the maps 
\begin{gather*}
\alpha = \begin{pmatrix} \alpha_{11} & 0 \\ 0 & \alpha_{22} \end{pmatrix}, \qquad
\beta = \begin{pmatrix} \beta_{11} & 0 \\ 0 & \beta_{22} \end{pmatrix}, \qquad
\tbeta = \begin{pmatrix} 0 & \beta_{21} \\ \beta_{12} & 0 \end{pmatrix}, \\
\gamma = \begin{pmatrix} \gamma_1 & \gamma_2 \end{pmatrix}, \qquad
\eta = \begin{pmatrix} \eta_1 \\ \eta_2 \end{pmatrix}
\end{gather*}
between the 1-cells $P=P_1\oplus P_2$ and $Q=Q_1\oplus Q_2$ in $C_V$.
Compatibility is ensured by the fact that the 2-cells of $M$ live inside $V$, and the vertical drop 
axioms therefore need not be verified separately. 

$M$ has a weak involution $\star$ interchanging $P$ and $Q$, and also transforming
$$
\alpha \mapsto -\alpha + \eta\circ\gamma, \qquad\qquad \beta \mapsto \tbeta, 
\qquad\qquad \gamma \mapsto \eta. 
$$
The 2-category $M$ is in fact fiat; indeed, there only remains to verify the existence of a bi-adjunction between $P$ and $Q$. 
One has a unit and co-unit, 
$$
\xi = 
\begin{pmatrix} 
\eta \\ \beta\circ\eta
\end{pmatrix} 
\colon I \to Q\oplus P 
\qquad\text{and}\qquad 
\epsilon = 
\begin{pmatrix}
\gamma & - \gamma\circ\beta
\end{pmatrix}
\colon P\oplus Q \to I,
$$
respectively. For instance, 
\begin{multline*}
\epsilon\iota_Q \circ \iota_Q\xi = 
\begin{pmatrix}
Q^*(\gamma) & - Q^*(\gamma)\circ Q^*(\beta)
\end{pmatrix}
\circ
\begin{pmatrix} 
Q_*(\eta) \\ Q_*(\beta)\circ Q_*(\eta)
\end{pmatrix} 
\\
=
\begin{pmatrix}
\iota_Q + \alpha & - (\iota_Q+\alpha)\circ (\alpha+\beta)
\end{pmatrix}
\circ
\begin{pmatrix} 
\iota_Q+\tilde\beta \\ \iota_P\circ (\iota_Q+\tilde\beta)
\end{pmatrix} 
\\
=
\begin{pmatrix}
\iota_Q+\alpha & - \alpha
\end{pmatrix}
\circ 
\begin{pmatrix} 
\iota_Q + \tilde\beta \\ \tilde\beta
\end{pmatrix} 
= \iota_Q,
\end{multline*}
and one similarly verifies $\iota_Q\epsilon \circ \xi\iota_Q = \iota_Q$ and 
$\epsilon\iota_P \circ \iota_P\xi = \iota_P = \iota_P\epsilon \circ \xi\iota_P$.

Now, the left cells of $M$ are $\{I\}$ and $\{P,Q\}$. 
The Duflo involution of this latter left cell clearly cannot be self-dual with respect to $\star$, since neither $P$ nor $Q$ is. 
We are appreciative of the assistance received from Dr Anne-Laure Thiel in simplifying this example. 
\eex

\end{document}